\documentclass[11pt,a4paper]{article}
\usepackage{amssymb,amsmath}
\usepackage{amsthm}
\usepackage{bbm}
\usepackage{stmaryrd}
\usepackage{hyperref}
\usepackage[usenames,dvipsnames]{xcolor}
\usepackage{color}
\usepackage{tikzsymbols}
\usepackage{lineno}

\input xy
\xyoption{all} 

\numberwithin{equation}{section}

\newtheorem{theorem}{Theorem}[section]

\theoremstyle{definition}
\newtheorem{definition}[theorem]{Definition}
\newtheorem{remark}[theorem]{Remark}
\newtheorem{example}[theorem]{Example}

\newcommand{\Id}{\mathbbmss{1}}

\newcommand{\InHom}{\mbox{$\underline{\Hom}$}}

\newcommand{\rmi}{ \textnormal{i}}
\newcommand{\rmd}{\textnormal{d}}

\newcommand{\opp}{\textnormal{op}}

\DeclareMathOperator{\Span}{Span}

\DeclareMathOperator{\Ber}{Ber}

\DeclareMathOperator{\Hom}{Hom}

\newcommand{\catname}[1]{\textnormal{\texttt{#1}}}

\font\black=cmbx10 \font\sblack=cmbx7 \font\ssblack=cmbx5 \font\blackital=cmmib10  \skewchar\blackital='177
\font\sblackital=cmmib7 \skewchar\sblackital='177 \font\ssblackital=cmmib5 \skewchar\ssblackital='177
\font\sanss=cmss10 \font\ssanss=cmss8 
\font\sssanss=cmss8 scaled 600 \font\blackboard=msbm10 \font\sblackboard=msbm7 \font\ssblackboard=msbm5
\font\caligr=eusm10 \font\scaligr=eusm7 \font\sscaligr=eusm5  \font\fraktur=eufm10
\font\sfraktur=eufm7 \font\ssfraktur=eufm5 
\font\bsymb=cmsy10 scaled\magstep2
\def\all#1{\setbox0=\hbox{\lower1.5pt\hbox{\bsymb
       \char"38}}\setbox1=\hbox{$_{#1}$} \box0\lower2pt\box1\;}
\def\exi#1{\setbox0=\hbox{\lower1.5pt\hbox{\bsymb \char"39}}
       \setbox1=\hbox{$_{#1}$} \box0\lower2pt\box1\;}

\def\tx#1{{\fam0\relax#1}}

\newfam\bifam
\textfont\bifam=\blackital \scriptfont\bifam=\sblackital \scriptscriptfont\bifam=\ssblackital

\newfam\blfam
\textfont\blfam=\black \scriptfont\blfam=\sblack \scriptscriptfont\blfam=\ssblack

\newfam\bbfam
\textfont\bbfam=\blackboard \scriptfont\bbfam=\sblackboard \scriptscriptfont\bbfam=\ssblackboard

\newfam\ssfam
\textfont\ssfam=\sanss \scriptfont\ssfam=\ssanss \scriptscriptfont\ssfam=\sssanss
\def\sss#1{{\fam\ssfam\relax#1}}

\newfam\clfam
\textfont\clfam=\caligr \scriptfont\clfam=\scaligr \scriptscriptfont\clfam=\sscaligr

\newfam\frfam
\textfont\frfam=\fraktur \scriptfont\frfam=\sfraktur \scriptscriptfont\frfam=\ssfraktur

\def\hpb#1{\setbox0=\hbox{${#1}$}
    \copy0 \kern-\wd0 \kern.2pt \box0}
\def\vpb#1{\setbox0=\hbox{${#1}$}
    \copy0 \kern-\wd0 \raise.08pt \box0}

\def\pmb#1{\setbox0\hbox{${#1}$} \copy0 \kern-\wd0 \kern.2pt \box0}
\def\pmbb#1{\setbox0\hbox{${#1}$} \copy0 \kern-\wd0
      \kern.2pt \copy0 \kern-\wd0 \kern.2pt \box0}
\def\pmbbb#1{\setbox0\hbox{${#1}$} \copy0 \kern-\wd0
      \kern.2pt \copy0 \kern-\wd0 \kern.2pt
    \copy0 \kern-\wd0 \kern.2pt \box0}
\def\pmxb#1{\setbox0\hbox{${#1}$} \copy0 \kern-\wd0
      \kern.2pt \copy0 \kern-\wd0 \kern.2pt
      \copy0 \kern-\wd0 \kern.2pt \copy0 \kern-\wd0 \kern.2pt \box0}
\def\pmxbb#1{\setbox0\hbox{${#1}$} \copy0 \kern-\wd0 \kern.2pt
      \copy0 \kern-\wd0 \kern.2pt
      \copy0 \kern-\wd0 \kern.2pt \copy0 \kern-\wd0 \kern.2pt
      \copy0 \kern-\wd0 \kern.2pt \box0}

\mathchardef\za="710B  
\mathchardef\zb="710C  
\mathchardef\zg="710D  
\mathchardef\zd="710E  
\mathchardef\zve="710F 
\mathchardef\zz="7110  
\mathchardef\zh="7111  
\mathchardef\zvy="7112 
\mathchardef\zi="7113  
\mathchardef\zk="7114  
\mathchardef\zl="7115  
\mathchardef\zm="7116  
\mathchardef\zn="7117  
\mathchardef\zx="7118  
\mathchardef\zp="7119  
\mathchardef\zr="711A  
\mathchardef\zs="711B  
\mathchardef\zt="711C  
\mathchardef\zu="711D  
\mathchardef\zvf="711E 
\mathchardef\zq="711F  
\mathchardef\zc="7120  
\mathchardef\zw="7121  
\mathchardef\ze="7122  
\mathchardef\zy="7123  
\mathchardef\zf="7124  
\mathchardef\zvr="7125 
\mathchardef\zvs="7126 
\mathchardef\zf="7127  
\mathchardef\zG="7000  
\mathchardef\zD="7001  
\mathchardef\zY="7002  
\mathchardef\zL="7003  
\mathchardef\zX="7004  
\mathchardef\zP="7005  
\mathchardef\zS="7006  
\mathchardef\zU="7007  
\mathchardef\zF="7008  
\mathchardef\zW="700A  
\mathchardef\zC="7009  

\newcommand{\be}{\begin{equation}}
\newcommand{\ee}{\end{equation}}

\newcommand{\bea}{\begin{eqnarray}}
\newcommand{\eea}{\end{eqnarray}}
\newcommand{\Z}{{\mathbb Z}}
\def\*{{\textstyle *}}
\newcommand{\R}{{\mathbb R}}

\newcommand{\C}{{\mathbb C}}

\newcommand{\s}{{\textstyle *}}

\def\Hom{\sss{Hom}}

\def\Sec{\sss{Sec}}
\def\Vect{\sss{Vect}}

\def\sT{{\sss T}}

\def\xi{\tx{i}}

\def\s*{{\scriptstyle *}}

\def\cO{\mathcal{O}}

\newcommand{\beas}{\begin{eqnarray*}}
\newcommand{\eeas}{\end{eqnarray*}}

\def\half{\frac{1}{2}}

\begin{document}
\bibliographystyle{plain}
\author{Andrew James Bruce}
\title{A First Look at Supersymmetry}
\date{\today}

\maketitle

\abstract{These are expanded notes for a short series of lectures, presented at the University of Luxembourg in 2017, giving an introduction to some of the ideas of supersymmetry and supergeometry. In particular, we start from some motivating facts in physics, pass to the theory of supermanifolds, then to spinors, ending up at super-Minkowski space-times. We examine some salient  mathematical issues with understanding supersymmetry in a classical setting and make no attempt to discuss phenomenologically interesting models. Moreover, the presentation is ``light'' in the sense that nothing is carefully proved.  The audience of the seminars ranged from Ph.D. students to more experienced postdocs in mathematics. The audience was assumed to have a working knowledge of differential geometry, elementary category theory, and basic ideas from classical \& quantum mechanics. No knowledge of supersymmetry and supergeometry is assumed.  }
\tableofcontents

\bigskip 

\noindent \textbf{Disclaimer} \emph{These notes contain no original material and the author has made no attempt at providing a comprehensive list of references.}

\section{Prerequisites}
Given that these lectures are aimed at graduate students and postdocs in mathematics, we assume basic knowledge of the following:
\begin{itemize}
\item  General Topology - topologies, open sets and topological spaces.
\item  Differential Geometry - manifolds, vector fields and Lie groups.
\item  Category Theory - categories and functors.
\item  Lie Theory - Lie algebras and their representations.
\item  Classical Mechanics - Lagrangians, Hamiltonians and Equations of Motion.
\item  Quantum Mechanics - Hilbert spaces and Hermitian operators.
\end{itemize}
No particular deep results from the above will be employed, but a general working knowledge will be necessary. We will not need much from Algebraic Geometry and we will not assume any knowledge of `supermathematics'.  We will not deal with Quantum Field Theory at all.

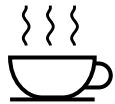
\begin{center}
\Coffeecup [5]
\end{center}

\section{Introduction}
In nature there appears to be two kinds of particles \footnote{There are other mathematical possibilities called \emph{Green-Volkov para-particles}, see \cite{Green:1953,Volkov:1959}. However, there is no evidence that (fundamental) particles with non-standard statistics are realised in nature.}: \emph{bosons} with integer spin, and \emph{fermions} with half-integer spin. These two families are very different when it comes to their statistics. In particular, any number of bosons can be placed into a given quantum state, while for fermions, we have the famous Pauli exclusion principle -- no more that one fermion can occupy a given quantum state.\par
At the classical level we know that bosons can be described by commuting objects, while fermions require anticommuting objects. The reason for the use of anticommuting objects can be quickly justified from the Palui exclusion principle and the classical limit of the canonical anticommutation relations. Throughout these notes, we will, generally, employ natural coordinates in which $c= \hbar =  1$. \par 
Consider a Hilbert space $\mathcal{H}_{F}$, which we equip with creation and annihilation operators $c^{\dag}$ and $c$. Naturally, these operators act on the vacuum (lowest energy state) as
\begin{align*}
& c^{\dag} | 0\rangle  = | 1\rangle \,,\\
& c| 0\rangle  = \mathbf{0}\,,\\
& c | 1\rangle  = | 0\rangle\,,
\end{align*}
where $|1\rangle$ is the one-particle state.  As we insist that no more than one fermion can be placed in a given state, it must be the case that
$$c^{\dag} | 1\rangle  = \mathbf{0}\,,$$
which loosely means that forcing another fermion into an occupied state destroys the state.  A general state is then a linear combination of these two states -- the empty state and one filled state. The direct consequences of this are that
\begin{align*}
& (cc^\dag  + c^\dag c)|0\rangle  = c| 1\rangle = |0\rangle \,,\\
& (cc^\dag  + c^\dag c)|1\rangle  = c^\dag| 0\rangle = |1\rangle  \,.
\end{align*}
Then, assuming we have normalised the vacuum state, we see that $(cc^\dag  + c^\dag c) = \mathbf{1}$. The classical limit  --  rather hands-wavingly  -- means that we are looking for objects that satisfy $\zx^1 \zx^2 +  \zx^2 \zx^1= 0 $. Thus, we encounter anticommuting objects -- this leads us to consider so-called \emph{Grassmann algebras}. Clearly, one cannot find real numbers, other than zero, that represent these anticommuting objects, Dirac referred to `anticommuting classical numbers'.\par 
The basic idea of \emph{supersymmetry} (c.f. \cite{Gervais:1971,Golfan:1971,Salam:1974,Wess:1974}) is to mix bosons and fermions. Quantum mechanically, consider a  Hilbert space with a decomposition into a bosonic and fermionic subspaces
$$\mathcal{H} = \mathcal{H}_{B} \oplus \mathcal{H}_{F}\,.$$ 
Now, suppose that we have self-adjoint operators  $Q_{I} = Q_{I}^{\dag}$, 
$(I = 1,2,\cdots , N)$, such that 
\begin{align*}
Q_{I}\mathcal{H}_{B} \subseteq \mathcal{H}_{F}, && Q_{I}\mathcal{H}_{F} \subseteq \mathcal{H}_{B}\,,
\end{align*}
and satisfy the \emph{supersymmetry algebra}
\begin{subequations}
\begin{align}
&\{Q_{I}, Q_{J} \} \propto \delta_{IJ} H\,,\\
& [Q_{I}, H] =0\,,\\
&[H,H] =0\,,
\end{align}
\end{subequations}
where $[A,B] = AB -BA$ is the \emph{commutator}, and   $\{ A, B\} = A B + B A$ is the \emph{anticommutator}. Here $H$ is the \emph{Hamiltonian} operator for our system under study. Some immediate remarks are:
\begin{enumerate}
\item The $Q$'s are the ``square-roots'' of the Hamiltonian. This means that supersymmetry is related to a space-time translation --  you apply one of the generators $Q$ twice and you end up acting with the Hamiltonian (the generator of time evolution of our system).
\item We can write $H = 2 Q^{\dag}_{I} Q_{I}$ (no sum here) and so the expectation value of the energy is
$$\langle  \psi | H | \psi \rangle = 2 \langle  \psi | Q^{\dag}_{I} Q_{I} | \psi \rangle  = 2 ||  Q_{I}|\psi \rangle  ||^2  \geq 0\,.$$
Thus, the energy spectrum is always non-negative. The vacuum energy is an order  parameter for symmetry breaking. That is, if the vacuum (ground) state is supersymmetric then
$$Q_{I}| 0 \rangle =\mathbf{0}\, , $$
and the vacuum energy is zero. Otherwise the energy of the vacuum is non-zero (but positive).
 \end{enumerate}
The simplest supersymmetric quantum mechanical theory one can think of is the $N=1$ \emph{supersymmetric oscillator} (see \cite{Nicolai:1976}). Here the Hamiltonian is simply the sum of the standard bosonic oscillator and the standard fermionic oscillator. In particular, on the Hilbert space we have bosonic creation and annihilation operators $a^\dag$ and $a$, and fermionic creation and annihilation operators $c^\dag$ and $c$ subject to the standard commutation/anticommutation relations
\begin{align*}
&[a, a^{\dag}] = \mathbf{1}, && \{ c, c^{\dag} \} =\mathbf{1}\,.
\end{align*} 
We supplement this with the condition $[a,c] =\mathbf{0}$, which basically says that the two sectors do not interact. The Hamiltonian is taken to be
$$H = a^\dag a + c^\dag c\,.$$
The SUSY generator is then given by
$$Q = c^\dag a + a^\dag c.$$
 One can show that as $a^\dag |0\rangle$ is a one boson state that $Qa^\dag | 0\rangle $ is a one fermion state and so on.\par 
 It must be stated that to date, supersymmetry has not been found in nature. However, methods of supersymmetry have been applied to condensed matter physics, nuclear and atomic physics, and optics. Wider than that, supersymmetry has inspired many developments in mathematics, especially in geometry as we shall see. \par 
 In practice writing down interacting supersymmetric quantum theories is difficult, if not impossible. A better tactic is to examine classical theories and then quantise these -- say via path integrals and perturbation theory. As the classical limit of a fermionic degree of freedom is an anticommuting object (think ``coordinate''), \emph{supermanifolds} become essential. To quote Witten (see \cite{Witten:2019})
 \begin{quote}
 Formulating a supersymmetric field theory in superspace – that is on a supermanifold – is, when possible, often very helpful.
 \end{quote}
Recent advances in using supermanifolds to build and understand supersymmetric theories include \cite{Castellani:2014,Castellani:2015,Castellani:2015b}. We remark that the literature is constantly growing. \par 
The idea of this short series of lectures is to support the idea that supermanifolds are useful in physics,  and to make sense of some of the mathematical gaps usually found in introductory physics texts. Very loosely, if we consider classical configuration or phase spaces of systems that include fermions, we require anticommuting coordinates and this leads to the idea of supermanifolds. Heuristically, supermanifolds are `manifolds' with commuting and anticommuting coordinates. We will further motivate the use of anticommuting coordinates and describe supersymmetry,  before passing to supermanifolds a bit more carefully. 

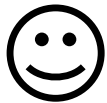
\begin{center}
\Smiley [5]
\end{center}
 
\section{Supersymmetric Mechanics}\label{Sec:mechanics}
\subsection{N=1 supermechanics in superspace}
As already stated, the idea is to find some classical action that is invariant under supersymmetry and then pass that over to the `black box' of quantisation. In these lectures I will say nothing about the quantum theory and just highlight some of the supergeometric aspects of the classical theory. In order to build some intuition, and find some interesting mathematical gaps, we first concentrate on a simple model in one dimension. This will allow us to highlight some issues to cover in later lectures without the clutter of indices and knowledge of relativity and spinors. In fact, all the subtle issues are already present in this simple model. \par 
Let us start from the supersymmetry algebra (see Appendix \ref{appendix:Lie} for Lie superalgebras)
\begin{align}\label{eqn:N=1}
\{ Q, Q\} = \frac{1}{2}P\,,&& [P,P] =0\,,  &&  [Q,P] =0\,.
\end{align}
Note that I have changed the normalisation here as compared with earlier.  This is rather unimportant and other conventions can be found in the literature.  The fundamental idea is to realise this algebra as vector fields on some supermanifold -- a concept we introduce properly in the next lecture. Indeed, the algebra here is a Lie superalgebra, so a  $\Z_2$-graded generalisation of a Lie algebra. Thus, one might expect to be able to integrate this to a `Lie supergroup'. This can be done, but requires some mathematics that we will discuss in a later lecture.  For now, as is common in physics, lets work formally and not worry about the exact mathematical meaning.\par 
I claim that one can realise $Q$ and $P$ as vector fields on the supermanifold $\R^{1|1}$. Without a proper definition, we can think of $\R^{1|1}$ as a `manifold' with one commuting and one anticommuting coordinate, $(t, \theta)$, say. In particular, we have $t \theta =  \theta t$ and  $\theta^2 = \theta \theta =  - \theta \theta = 0$. And so, all functions on $\R^{1|1}$ are simply of the form
$$f(t,\theta) = f_{0}(t) + \theta f_{1}(t)\,.$$ 
One can justify this from several points of view, but we can think of Taylor expanding a function in $\theta$ and noting the nilpotent property of anticommuting coordinates, i.e., $\theta^2 =0$ and so higher-order terms in any expansion vanish.  Note that a general function decomposes into a bosonic and fermonic part.  If we assign the number $0$ to anything bosonic and the number $1$ to anything fermionic then we observe that the space of functions is $\Z_2$-graded\footnote{Recall that as a set $\Z_2 := \{ 0,1\}$. }. As an aside, by supergeometry in a general sense, we mean geometry in which everything is $\Z_2$-graded --  this is not always an obvious extension of geometry!     \par 
Now, one can show that
\begin{align}
P = \frac{\partial}{\partial t}, && Q = \frac{\partial}{\partial \theta} + \frac{1}{4} \theta \frac{\partial}{\partial t},
\end{align}
gives a representation of (\ref{eqn:N=1}), i.e., these vector fields (first-order differential operators) satisfy the supersymmetry algebra. Here we take the derivative with respect to $t$ to have it's standard meaning, while derivatives with respect to $\theta$  are defined algebraically -- functions are just linear in $\theta$ and so this is no problem. In more generallity, functions in odd/anticommuting variables will be polynomials, and so derivatives can be defined algebraically.\par
We need to be a little bit careful here with the commutator of vector fields as we are in the `super-world' -  everything is $\Z_2$-graded. We can deduce the parity, so if a vector field is even/bosonic or odd/fermionic, by counting the number of anticommuting coordinates and their derivatives (mod 2 when required). Thus, $P$ is even/bosonic and $Q$ is odd/fermionic. This means that we should consider the standard commutator between $P$ and itself, and between  $P$ and  $Q$. However, we need the anticommutator between $Q$ and itself, i.e.
$$\{Q,Q\} = QQ +QQ\,. $$  \par 
As we have a vector field that represents the supersymmetry generator, we can think about its Lie derivative.  In particular, allowing an anticommuting infinitesimal parameter $\epsilon$, we can geometrically understand supersymmetry as  the infinitesimal diffeomorphism 
\begin{subequations}
\begin{align}
& t \mapsto t' = t +\frac{1}{4} \epsilon \theta\,,\\
& \theta \mapsto \theta' = \theta + \epsilon\,.
\end{align}
\end{subequations}
But what has this got to do with our physics understanding of supersymmetry? The idea now is to encode the bosonic and fermionic degrees of freedom into one object on $\R^{1|1}$ -- a so called \emph{superfield}. The theory we will build will have one bosonic field $q(t)$ and one fermionic field $\psi(t)$. With this in mind let us consider an even superfield on $\R^{1|1}$, which in coordinates must be of the form 
\begin{equation} \label{eqn:Superfield1}
\Phi(t, \theta) =  q(t) + \theta \psi(t)\,.
\end{equation}
Note that we require this whole object to be even/commuting. This means that we require $\psi(t)$ to be anticommuting. This is of course natural in light of the spin-statistics theorem -- ``fermions have anticommuting fields''. However, this means that we \emph{cannot} think of $\Phi$ simply as a function on the superspace $\R^{1|1}$. How can we build odd/anticommuting fields that only depend on a commuting coordinate?  We can also consider odd superfields, but this does not remedy the situation. \par 
The resolution to this question of odd fields, as we have already done with understanding infinitesimal diffeomorphisms of $\R^{1|1}$, is to allow `external odd constants'. The best and most mathematically sound way to do this is to employ a more categorical framework to supermanifolds -- namely the functor of points and internal Homs. We will come to this in the next lecture. For now, let us agree that we can work formally and not worry about the extra oddness.\par 
In \emph{component form} we define the N=1 supersymmetry transformations via
$$\delta_{\epsilon}\Phi(t,\theta) :=  \epsilon Q(\Phi(t, \theta))\,,$$
which is really just the Lie derivative of the superfield in the `direction' of $Q$. Written out we have
\begin{subequations}
\begin{align}\label{eqn:N=1Comps}
& \delta_{\epsilon}q(t) = \epsilon \psi(t)\,,\\
\label{eqn:N=1Compsb}
 & \delta_{\epsilon}\psi(t) = - \frac{1}{4} \epsilon \dot q(t)\,,
\end{align}
\end{subequations}
where `dot' is the derivative with respect to $t$. Thus, we see that supersymmetry mixes the bosonic and fermionic degrees of freedom, as we wanted!\par 
We want to use this encoding of our fields/degrees of freedom into a geometric object on $\R^{1|1}$ to build a classical action that is invariant under supersymmetry -- that is, invariant under \eqref{eqn:N=1Comps} and  \eqref{eqn:N=1Compsb}. As we what something dynamical we should take derivatives of the superfield and build some action from there. \par 
The slight complication here is that derivatives with respect to $\theta$ transform non-trivially under supersymmetry -- check this yourself!  The situation is formally the same as encountered when introducing connections and covariant derivatives in differential geometry. Thus, one needs to introduce a \emph{SUSY covariant derivative}
\begin{equation}
\mathbb{D} :=  \frac{\partial}{\partial \theta} - \frac{1}{4} \theta \frac{\partial}{\partial t} \,.
\end{equation}
You should notice that this is almost the same as $Q$, but with a minus sign difference -- this is quite a generic feature. We want to build actions from the superfield and its first derivatives  (we will not consider higher-order derivatives as these are generally problematic in physics)
$$\left \{ \Phi, \dot \Phi, \mathbb{D}\Phi  \right\}\,.$$
Just about the only sensible action we can build is
$$S[\Phi] =  {-} \int D[t, \theta] \mathbb{D}\Phi \dot \Phi\,.$$
For odd variables integration is the same as differentiation (see Appendix \ref{App:calculus}), and then integrating out $\theta$ we are left with an action 
\begin{equation}
S = \int \rmd t \left( \frac{1}{4} \dot q^2(t) - \dot \psi(t) \psi(t) \right)\,,
\end{equation}
which gives (up to some choice in the mass) the standard kinetic term for a particle in one dimension and the one dimensional Dirac-like term. \par 
Via direct calculation you can show that the Lagrangian is quasi-invariant under supersymmetry:   $\delta_{\epsilon}L = \frac{1}{4} \epsilon \frac{\rmd}{\rmd t} \left( \dot q \psi\right)$. However, by construction we already knew that the action is invariant!
\begin{remark}
The one issue we have `brushed under the carpet' is the invariance or not of the integration measure under supersymmetry. To do this carefully, we need to introduce the \emph{Berezinian} which is the super-generalisation of the determinant.  At this stage, we will just state that the integration measure is invariant. For completeness, we include the calculation in Appendix \ref{appendix:Ber}.   
\end{remark}

\subsection{N=2 supermechanics in superspace}
The N=1 action we obtained is not particularly interesting -- it is just the sum of the free actions for a boson and a fermion. The N=2 case follows from similar reasoning, but we have a much more interesting theory.  I will just sketch the theory, making note that the basic mathematical questions already pose themselves for the N=1 case. \par 
The N=2 supersymmetry algebra is given by (noting our conventions) as
\begin{subequations}
\begin{align}
& \{ Q_{I}, Q_{I} \} = \frac{1}{2}\delta_{IJ} P\,,\\
& [Q_{I}, P] =0\,,\\
& [P,P]=0\,,
\end{align}
\end{subequations}
where $I = 1,2$. In this case we can take the supermanifold needed to realise this algebra to be $\R^{1|2}$. That is the supermanifold which we can equip with coordinates $(t, \theta^1, \theta^2)$. The $\theta$ coordinates are odd/anticommuting and so $\theta^1 \theta^2 = - \theta^1 \theta^2$ and $(\theta)^2 =0$. I claim, and one  should check (I mean you!) that we have the representation
\begin{align}
P = \frac{\partial}{\partial t}, && Q_{I} =  \frac{\partial}{\partial \theta} + \frac{1}{4} \theta^{J}\delta_{JI} \frac{\partial}{\partial t}\,.
\end{align}
Geometrically, we  understand supersymmetry as  the infinitesimal diffeomorphisms
\begin{subequations}
\begin{align}
& t \mapsto t' = t +\frac{1}{4} \epsilon^{I} \theta^{J}\delta_{JI},\\
& \theta^{I} \mapsto \theta'^{I} = \theta^{I} + \epsilon^{I}\,.
\end{align}
\end{subequations}
This means that we have two supersymmetries and require two anticommuting parameters. Again, just as for the N=1 case we need SUSY covariant derivatives
$$\mathbb{D}_{I} =  \frac{\partial}{\partial \theta^{I}} - \frac{1}{4} \theta^{J}\delta_{JI} \frac{\partial}{\partial t}\,.$$
An even superfield is then of the form
$$\Phi(t,\theta) =  q(t) + \theta^{1}\psi_{1}(t) +\theta^{2}\psi_{2}(t) + \theta^1 \theta^2 F(t)\,.$$
This we see that we have two bosonic degrees of freedom and two fermionic ones, as we would expect. It will turn out that $F$ is an auxiliary field, meaning it can be removed using its equations of motion. In component form the (off-shell) supersymmetry transformations are
\begin{subequations}
\begin{align}
\delta q(t) & = \epsilon^1 \psi_{1}(t) + \epsilon^2 \psi_{2}(t)\,,\\
\delta \psi_1(t) & = - \frac{1}{4} \epsilon^1 \dot q(t) + \epsilon^2 F(t)\,,\\
\delta \psi_2(t) & = - \frac{1}{4} \epsilon^2 \dot q(t) - \epsilon^1 F(t)\,,\\
\delta F(t) & = \frac{1}{4} \epsilon^1 \dot\psi_{2}(t) - \frac{1}{4} \epsilon^2 \dot\psi_{1}(t)\,.
 \end{align}
 \end{subequations}
 To build and invariant action we consider
 $$S[\Phi] =  \int D[t, \theta]  \left(\mathbb{D}_{1}\Phi \mathbb{D}_{2}\Phi  - W(\Phi) \right)\,.$$
 Expanding out in terms of the $\theta$'s and then integrating (see Appendix     \ref{App:calculus}) gives
 $$ S = \int \rmd t \left(\frac{1}{8} \dot q^2 - \frac{1}{4}\left(\dot \psi_{1} \psi_1 + \dot\psi_2 \psi_2 \right) + F^2 - F W^{\prime}(q) + \psi_1 \psi_2 W^{\prime \prime}(q) \right)\,.$$
 One can then remove $F $ using its equation of motion giving 
 \begin{equation}
 S =  \int \rmd t \left(\frac{1}{8} \dot q^2 - \frac{1}{4}\left(\dot \psi_{1} \psi_1 + \dot\psi_2 \psi_2 \right)  -  \frac{1}{4} \left(W^{\prime}(q) \right)^2  + \psi_1 \psi_2 W^{\prime \prime}(q)\right)\,,
 \end{equation}
 which is related to Witten's model of supersymmetric quantum mechanics (see \cite{Witten:1981}).

\subsection{Remarks}
I have been very sketchy in this first lecture. The key message is that supergeometry can help with the construction of classical actions that exhibit supersymmetry.  \par 
The main issues that we saw are the following:
\begin{itemize}
\item What are supermanifolds?
\item How do we make sense of superfields and especially odd fields?
\end{itemize}
I hope that the next lecture will answer these questions -- to some reasonable extent.

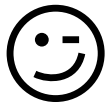
\begin{center}
\Winkey [5]
\end{center}

\section{Supermanifolds and the functor of points}
\subsection{Recap}
Supersymmetry at the quantum level relates bosons (integer spin particles) with fermions (half integer spin). As these particles are very different in the essence, the theoretical existence of such a symmetry is quite remarkable.   However, constructing interesting quantum theories from scratch is difficult and one usually resorts to classical theories first. The ``particle-field correspondence'' (together with the ``spin-statistics theorem'')  tells us
\begin{itemize}
\item Bosons  $\longrightarrow  q(t,x)$  \hspace{50pt} commuting fields
\item Fermions $\longrightarrow  \psi(t,x)$  \hspace{40pt} anticommuting fields
\end{itemize}
The questions we shall address in this lecture are the following
\begin{enumerate}
\item Can we \emph{geometrically} treat commuting and anticommuting degrees of freedom on equal footing?
\item Can we really make sense of $\psi(t,x)$ being anticommuting if it seems to be `build' only from \emph{commuting} variables?
\end{enumerate}
Motivation for these questions comes from classical mechanics and field theory where differential geometry plays a central r\^ole. We need some kind of `manifold' that has both commuting and anticommuting coordinates. Mathematicians took (hijacked?) the nomenclature \emph{super} from supersymmety and developed the theory of \emph{supergeometry}, which is now a respectable, if not a very popular, branch of mathematics. 

\subsection{Supermanifolds}
The very informal working definition of a supermanifold is a `manifold' with commuting and anticommuting coordinates. One can work with coordinate charts and patch them together, much as one would with a manifold. From there various geometric objects can be defined via their local description in terms of coordinates. Everything looks very much like the theory of manifolds, at first glance.  To make these loose statements more meaningful mathematicians developed at least two approaches to this. First is the DeWitt--Rogers approach where one gives a non-Hausdorff topology to Grassmann algebras and then builds supermanifolds from the according local model. There are some complications here with the general theory and in particular, one has to take some care in defining (different classes of) smooth functions and vector fields.  However, for application in physics, the DeWitt--Rogers approach, also known as the concrete approach, is seemingly sufficient. For a review the reader can consult Rogers \cite{Rogers:2007}.\par 
 The other is the ringed space approach were one follows a more algebraic geometry route and modifies the sheaf of (continuous) functions on a topological space. It is this ringed space approach that we will adopt as this is known to be more `mathematically sound' and is the favoured approach by mathematicians today. In fact, it is known that sorting out the issues of the concrete approach leads one back to sheaves, and so the initial advantage of a less abstract language is actually not there.  The ringed space approach goes back to Berezin (1931-1980) who is recognised as the father of supermathematics. While there had been some hints in physics literature prior to Berezin's contributions, he was the first to have the foresight to envisage treating commuting and anticommuting variables on equal footing.  The definition of a supermanifold, as we understand it, is due to Berezin and Le\u{\i}tes (see \cite{Berezin:1975,Leites:1980}). Other useful references include \cite{Carmeli:2011,Deligne:1999,Freed:1999,Manin:1997,Varadrajan:2004}.
\begin{definition}
 A (smooth) \emph{supermanifold} of dimension $n|m$ is a locally superringed space $M := \big( |M|, \cO_{M} \big)$, which is locally isomorphic to $\R^{n|m}:=  \big(  \R^{n}, C^{\infty}_{\R^{n}}\otimes \Lambda(\zx^{1}, \cdots , \zx^{m})\big)$, where  $C^{\infty}_{\R^{n}}$ is the sheaf of smooth functions on $\R^n$, and $\Lambda(\zx^{1}, \cdots , \zx^{m})$ is the Grassmann algebra (over the reals) with $m$ generators, i.e., the generators satisfy $\zx^i \zx^j = {-}\zx^j \zx^i$, which implies $(\zx)^2 =0$.  \emph{Morphisms of supermanifolds} are morphisms of locally ringed spaces,
 $$(\phi, \phi^{*}):  \big(  |M|, \cO_{M} \big) \longrightarrow \big(  |N|, \cO_{N} \big)$$
 Thus, they consist of two pieces:
 \begin{enumerate}
 \item A  continuous map $\phi: |M| \rightarrow |N|$;
 \item A family of superring morphism $ \phi^{*}_{|V|}: \cO_{N}(|V|)  \rightarrow \cO_{M}(\phi^{-1}(|V|))$, where $|V| \subset |N|$ is open, that commute with the restriction maps. Via composition of morphisms of supermanifolds we obtain the \emph{category of supermanifolds} $\catname{SMan}$
 \end{enumerate}
\end{definition}

\medskip 
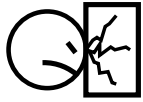
\begin{center}
\Walley [5]
\end{center}

\noindent Let us decipher some of this...\par 
\medskip 

In the definition of a supermanifold we have a ringed space of some kind, i.e., $|M|$ is a topological space. In fact, from the definition it turns out that $|M|$ is a \emph{smooth manifold}. One can justify this by `throwing away' the coordinates $\zx$ and observing that we have the definition of a manifold!  \par 
Given any smooth manifold, $|M|$, we can consider its sheaf of smooth functions. That is, given any open $|U| \subset |M|$,
$$\cO(|U|) = C^{\infty}(|U|, \R) =: C^{\infty}(|U|) \,.$$
 As any `small enough' open $|U|$ is homomorphic to some subset in $\R^{n}$, we really are considering multi-variable smooth functions, i.e., we can use local coordinates. Mod details, this means that we can always think of functions locally in a consistent way. The interested reader should look up the properties of sheaves, however we will not need any details here.\par 
The idea of a supermanifold is to modify the sheaf of smooth functions on a manifold to obtain a new kind of geometry -- one `throws in' generators of a Grassmann algebra
\begin{equation}\label{eqn:StructureSheaf}
\cO_{M}(|U|) \simeq  C^{\infty}(|U|) \otimes \Lambda(\zx^{1}, \cdots , \zx^{m})\,, 
\end{equation}
for `small enough' opens $|U|\subset |M|$. This means that local sections of the \emph{structure sheaf} of a supermanifold are of the form (using the isomorphism above and local coordinates on $|U|$)
$$f(x, \zx) = f(x) + \zx^{i}f_{i}(x) + \frac{1}{2!}\zx^{i}\zx^{j}f_{ji}(x) + \cdots \frac{1}{m!} \zx^{i_{1}} \cdots \zx^{i_{m}} f_{i_{m} \cdots i_{1}}(x).$$
We will write $C^{\infty}(M)$  for the global sections and call the \emph{functions}. Note that any function written in local coordinates is a polynomial of maximal order $m$.  The dimension of a supermanifold is $n|m$, where $n$ is the dimension of the underlying manifold, and $m$ is the number of generators of the Grassmann algebra.  \par 
As a warning, on some arbitrary open $|U| \subset |M|$ we know that $\cO_M(|U |)$ is a (supercommutative) unital superring, but that that is about it. It is only on  `small enough' opens can we use coordinates - this is exactly the same as with manifolds.\par 
Another thing to notice is that we need two parts to every morphism of supermanifolds. For smooth manifolds one can show that the only possible pullback of functions that is meaningful is that induced by the standard pullback by smooth maps between the manifolds. However, for supermanifolds, as we have modified the notion of functions this is no longer the case and the pullback needs to be specified. \par 
What makes supermanifolds very workable from a differential geometry and physics perspective is that we can use \emph{local coordinates} in much the same way as we can on a smooth manifold. Indeed, much of the theory of smooth manifolds can be directly translated to supermanifolds via local coordinates -- not that this is always very elegant. 
\begin{theorem}[``Chart Theorem'']
We can use local coordinates on a supermanifold and morphisms of supermanifolds can be described locally by their effect on coordinates.
\end{theorem}
That is, given any point on $|M|$ we can always find a `small enough' open neighbourhood $|U|\subseteq |M|$ such that we can  employ local coordinates $x^{a} := (x^{\mu} , \theta^{i})$ on $M$, where $x^{\mu}$  and $\theta^{i}$ are, respectively, collections of commuting and  anticommuting elements of $\cO_{M}(|U|)$. Using the isomorphism (\ref{eqn:StructureSheaf}) on can set $x^{\mu}$ to be a `standard coordinate' on $|U| \subset \R^{n}$ and $\theta^i$ are the generators of a Grassmann algebra. However, it is not usually instructive to think of $x$ as `pure real' and $\theta$ as `pure Grassmann' as we have changes of coordinates that mix these two, as we will see by looking at more general morphisms.\par 
Now suppose we have another supermanifold $N$ and we equip it with local coordinates $(y^{\lambda} , \chi^{\delta})$, the first being commuting coordinates and the latter being anticommuting coordinates. Then any morphism, which is usually written simply as
$$\phi : M \rightarrow N\,,$$
can be locally described as
\begin{align*}
& \phi^{*}y^{\lambda} = y^{\lambda}(x, \theta)\,,\\
& \phi^{*} \chi^{\delta} = \chi^{\delta}(x, \theta)\,,
\end{align*}
where we have employed standard abuses of notation form standard differential geometry. Moreover, from the definition of a supermanifold we require that the morphism must preserve the \emph{Grassmann parity}, i.e., the commuting or anticommuting nature of the local coordinates. 
\begin{example}
Let us consider $\R^{1|2} \longrightarrow \R^{1|1}$. In this case we can consider global coordinates $(x, \theta^{1} , \theta^{2})$ on $\R^{1|2}$, and $(t, \tau)$ on $\R^{1|1}$. Then using the standard abused of notation we write
\begin{align*}
& \phi^{*}t = t(x, \theta) = t(x) + \theta^{1}\theta^{2}s(x)\,,\\
& \phi^{*}\tau = \tau(x, \theta) = \theta^{1}\tau_{1}(x) + \theta^{2}\tau_{2}(x)\,,
\end{align*}
where $t(x), s(x), \tau_{1}(x)$ and $\tau_{2}(x)$ are smooth functions on the real line.
\end{example}
We remark that the ``Chart Theorem'' is a non-trivial statement. The fact that we have local coordinates follows directly from the definition of a supermanifold. The fact that morphisms of algebras are determined by morphisms of their generators in \emph{part} establishes the theorem. In particular, if we were only dealing with polynomials in both even and odd coordinates then there theorem would be rather obvious: knowing what happens to the coordinates means we know what happens to the polynomials, and vice versa. However, we are allowing more general smooth functions in the commuting coordinates and this needs some subtle handling in order to establish the theorem -- one uses the  polynomial approximation to prove the theorem.  We will not attempt prove in any rigorous way the ``Chart Theorem'' and simply take it as given.\par 
The upshot of the local structure and chart theorem is that you can think of a supermanifold as a space with both commuting and anticommuting coordinates...

\medskip
 
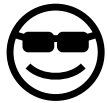
\begin{center}
\Cooley [5]
\end{center}

 \medskip 
 
 ... almost.

\subsection{The functor of points}
The only genuine points on a supermanifold $M$ are the points on the underlying (or reduced) manifold $|M|$. Thus, not all the information about the supermanifold is contained in the set of point of $|M|$. This is in stark contrast with the classical case of manifolds. Loosely, once you know all the maps $* \mapsto |M|$, you know the manifold. However, one has to be careful with point-wise constructions in the category of supermanifolds. However, we can clawback some classical intuition using Grothendiech's \emph{functor of points}.\par 
The idea follows from more categorical nonsense (see Maclane \cite{MacLane1998})  and the \emph{Yoneda embedding}, which works for all small categories (i.e., the hom-sets are sets). One considers a supermanifold to be a functor in the functor category  $\widehat{\catname{SMan}} :=  [\catname{SMan}^\opp , \catname{Set}]$, i.e.,
\begin{align*}
M : & \, \catname{SMan}^\opp \longrightarrow \catname{Set}\\
    & S  \mapsto  M(S) := \Hom(S,M)
\end{align*} 
and  for any $\phi \in \Hom(P,S)$ we have
\begin{align*}
M(\phi) : &  M(S) \longrightarrow M(P)\\
    & m  \mapsto  m' := m \circ \phi.
\end{align*}
\emph{Yoneda's lemma} says that this assignment of a functor to a supermanifold is full and faithful,  hence $\catname{SMan}$ is a full subcategory inside $\widehat{\catname{SMan}}$. In particular, morphisms of supermanifolds correspond to natural transformations between the associated functors -- which really then reduce to considering maps between Hom sets.\par
 Basically, no information is lost by thinking of a supermanifold as a functor. In fact, in general working with functor categories allows for ideas and constructions to be employed that are just not available in the original category. The big advantage with thinking of supermanifolds as functors is that we are now dealing with sets and this can allow us to start to work a bit more classically in terms of `generalised points on a suprmanifold' known as  \emph{S-points}. However, in general, some care must be taken here as $M(S)$ for a given supermanifold $S$ is just a set and has no further structure. It is only by full knowledge of $M(-)$ as a functor, so how it acts on all supermanifolds (and morphisms) can we really fully recover $M$.
 \begin{remark}There is a nice particle physics analogy here with the functor of points, which we think is originally due to Ravi Vakil, though I stand to be corrected. Imagine that we want to discover all the properties of an unknown particle, and to this end we fire at it \emph{all} possible test particles with a wide range of energies.  By observing the results of the  interactions   with \emph{all} the test particles over \emph{all} energy scales, we can deduce \emph{all} the properties of the unknown particle and thus characterise it completely. The functor of points allows us to `probe' a given supermanifold by observing how, as a functor, it `interacts' with all `test supermanifolds'. We further remark that the nomenclature `probe' in this context goes back (at least) to Urs Schreiber pre-2008 online. 
\end{remark}
Clearly, not all functors in $\widehat{\catname{SMan}}$ come from the Yoneda embedding. We say that a functor in  $\widehat{\catname{SMan}}$ is \emph{representable} if it is naturally isomorphic to a supermanifold in the image of the Yoneda embedding -- in simpler terms it is a genuine supermanifold. We will then refer to general (i.e., not necessarily representable) functors in $\widehat{\catname{SMan}}$ as \emph{generalised supermanifolds}. In practice, it maybe easier to define some construction of a supermanifold as a generalised supermanifold and later show representability -- this is also often the case in algebraic geometry more generally.  \par 

\subsection{The Cartesian product}
The Cartesian product of supermanifolds can be constructed using local coordinates in more-or-less the same way as one does for manifolds. We will just sketch the construction here. \par

Let $M = (|M|, \cO_M)$ and $N = (|N|,\cO_N)$ be two supermanifolds of dimension $p|q$ and $r|s$, respectively. First note that the Cartesian products $|U|\times |V|$, where $|U|\in |M|$, $|V|\in |N|$ are open, form a basis of the product topology of $|M|\times |N|$. We then pick coverings of $M$ and $N$, labelling  the opens as $|U_i|$ and $|V_i|$, such that we can employ local coordinates $(x_i, \zx_i)$ and $(y_j, \eta_j)$ on $M$ and $N$, respectively. We denote the basis of the product topology generated by $|U_i|\times |V_i|$  as $\mathcal{B}$.  
\begin{definition}
Let $M$ and $N$ be supermanifolds of dimension $p|q$ and $r|s$, respectively. The \emph{product supermanifold} $M\times N$, of dimension $p+r| q+s$, is the supermanifold $(|M|\times |N|, \cO_{M\times N})$, where $|M|\times |N|$ is equipped with the product topology, and $\cO_{M\times N}$ is the sheaf glued from the sheaves $C^\infty_{|U_i|\times|V_j|}(x_i , y_j)\otimes\Lambda(\zx_i, \eta_j)$ associated to the base $\mathcal{B}$, i.e.,
$$\cO_{M\times N}|_{|U_i|\times  |V_j|} \simeq C^\infty_{|U_i|\times|V_j|}(x_i , y_j)\otimes\Lambda(\zx_i, \eta_j)\,.$$
\end{definition} 
The upshot of the construction is that we can employ local coordinates on $M$ and $N$ to build the Cartesian product $M\times N$. It can be shown that the Cartesian product is a categorical product (for those that know what that means). \par 
An important and useful result is the canonical identification of the functor of points
$$(M\times N)(S) = M(S) \times N(S)\,.$$
One application of the functor of points and the Cartesian product are \emph{super Lie groups}, which can be defined a group objects in the category of supermanifolds. To define this carefully, which we won't here, required the Cartesian product.  Equivalently, we can define a super Lie group as a functor
$$G : \catname{SMan}^{\opp} \longrightarrow \catname{Grp}\,,$$
such that composition  with the forgetful functor
$$\catname{Grp} \longrightarrow \catname{Set}\,,$$
is a representable functor,  that is, $G$ is a genuine supermanifold. In simpler language, every $G(S)$ is a group and changing parametrisation yields and group homomorphism.  We will encounter super Lie groups when dealing with supersymmetry.\par 
We also require the Cartesian product for generalised supermanifolds and enriched morphisms of supermanifolds. 

\subsection{Generalised supermanifolds and the inner Homs}
 As already stated, a generalised supermanifold is an object in the functor category $\widehat{\catname{SMan}}$, whose objects are functors from $\catname{SMan}^{\opp}$ to the category $\catname{Set}$ and whose morphisms are natural transformations.  One passes from the category of finite--dimensional supermanifolds to the larger category of generalised supermanifolds in order to understand, for example, the \emph{internal Homs} objects. In particular there always exists a generalised supermanifold such that the so--called \emph{adjunction formula} holds
$$\InHom(M,N)(\bullet) := \Hom(\bullet \times M,N) \in  \widehat{\catname{SMan}}\, .$$
Heuristically, one should think of enriching the morphisms between supermanifolds to now have the structure of a supermanifold, however to understand this one passes to a larger category. In essence we will use the above to define what we mean by a \emph{mapping supermanifold} and will probe it using the functor of points. We will refer to `elements' of a mapping supermanfold as \emph{maps} reserving morphisms for the categorical morphisms of supermanifolds.\par
A little more carefully, any arbitrary supermanifold $S$ is sent to the set
$$\InHom(M,N)(S) = \Hom(S \times M,N)\,,$$
and any given  morphism $\phi \in \Hom(P, S)$ induces the mapping between sets
\begin{eqnarray*}
 \InHom(M,N)(\phi) :   \InHom(M,N)(S)& \longrightarrow& \InHom(M,N)(P)\, ,\\
 \Psi_{S} &\longmapsto&  \Psi'_{P}:= \Psi_{S}\circ (\phi, \Id_{M}) \, .
\end{eqnarray*}
In general, the mapping supermanifold will not be representable. \par
We can think of maps between two supermanifolds as a family of arrows $M \longrightarrow N$ consisting of both even and odd maps, or as standard supermanifold morphisms $S\times M \longrightarrow N$ parametrised by arbitrary supermanifolds $S$.\par 
 To clarify this, let us employ local coordinates $(x^\mu , \theta^{i})$ on $M$, $(y^{\lambda}, \zx^{\delta})$ on $N$ and $(s^{t}, \chi^{v})$ on $S$ (chosen arbitrarily). Then, any map $\Phi_{S} \in \Hom(S\times M, N)$ is locally   (using standard abuses of notation) described by
$$  y^{\lambda} = y^{\lambda}(z, \chi, x, \theta) =:  y^{\lambda}_{S}(x, \theta)\, , \quad \zx^{\delta} =  \zx^{\delta}(s, \chi , x , \theta) =: \zx^{\delta}_{S}(x,\theta )\, ,$$ 
which are understood as collections of (local) even and odd functions on the Cartesian product $S\times M$, respectively. Note that with respect to the local coordinates on $M$ the maps maybe either even or odd: it is the Grassmann parity of the local coordinates on $S$  that `corrects' this discrepancy. Thus, the idea of even and odd maps is realised.
\begin{example}
 It is easy to see that $\Hom(\{\emptyset\},M) = |M|$, while $\InHom(\{\emptyset\},M) = M$.
\end{example} 
\begin{example}\label{exAntTangBund}
  Another well--known example of a representable generalised supermanifold is the antitangent bundle $\InHom(\mathbb{R}^{0|1},M) = \Pi \sT M$. In fact, it is well--known that the generalised supermanifolds $\InHom(\mathbb{R}^{0|p}, M)$  are representable for all $p \in \mathbb{N}$.
\end{example}
As maps between supermanifolds are really functors, the composition of maps should be thought of in terms of a natural transformation
\begin{equation}
\underline{\circ} : \InHom(M,N) \times \InHom(N,L) \longrightarrow \InHom(M,L)\, ,
\end{equation}
defined, for any $S \in \catname{SMan}$,  by
\begin{eqnarray}
 \Hom(S\times M, N) \times \Hom(S\times N,L)& \longrightarrow& \Hom(S\times M, L)\, \\
  (\Phi_{S}, \Psi_{S})& \longmapsto &(\Psi \,\underline{\circ} \,\Phi)_{S}:= \Psi_{S} \circ (\Id_{S}, \Phi_{S})\circ (\Delta, \Id_{M})\, ,\nonumber
\end{eqnarray}
and then letting $S$ vary   over \emph{all} supermanifolds. Here $\Delta:S\longrightarrow S\times S$ is the diagonal of $S$ and $\Id_S:S\longrightarrow S$ is the identity.\par 

\subsection{Anticommuting fields revisited}
Let us again reconsider mechanics with  anticommuting degrees of freedom. We would like to make sense of anticommuting objects that depend on time only. In general, we will want to understand fields that depend on space also. For now, we want to understand how $\psi^{\alpha}(t)$ can be made anticommuting (see for example (\ref{eqn:Superfield1})). Clearly we cannot understand this as a categorical map $\psi : \R \rightarrow \R^{0|m}$ as the pullback of the odd coordinates on $\R^{0|m}$ to functions on $\R$ will vanish -- the only numerical value the generators of a Grassmann algebra can take are zero. We thus need an \emph{odd map}.\par 
One solution could be to allow anticommuting coefficients in the theory. That is, we will allow `constants' that take their values in some `external' Grassmann algebra. But, what Grassmann algebra should we pick? If we find one that is `large enough' to describe the theory, then why not consider a larger one? If we expand $\psi^{\alpha}(t)$ in terms of these chosen parameters, then do the components have any physical meaning? \par 
I argue that the choice of odd constants has no physical meaning -- they are there just to ensure that everything is of the right nature with respect to being commutative or anticommutative. None of the physics should depend in any critical way on the choice of odd constants. Thus, don't make a choice! \par 
That is, we anticommuting (fermion) fields as \emph{parametrised functors},  in particular in our mechanical situation they are functors from $\catname{SM}^{\opp}$  to $\catname{Set}$ that are parametrised by time: so we use the internal homs and write (rather abusing set theory notation) ``$\Psi \in \InHom(\R, \R^{0|m})$''. To make sense of this one `probes'  the generalised supermanifold of fields and thinks of
$$\Psi_{S} : S \times \R \longrightarrow \R^{0|m}\,.$$
Let us then equip $S$ with local coordinates $(s^{t}, \chi^{v})$, then one can 
 write (via pulling back the coordinates on $\R^{0|m}$)
$$\psi^{\alpha}(s, \chi, t) := \psi^{\alpha}_{S}(t).$$
Thus the extra `anticommutingness' is provided for by the local coordinates on $S$. In general, we don't usually write out the coordinates on $S$ explicitly, but it can be useful to do so when keeping track of some computations -- in fact we have a theorem by Schwarz and Voronov \cite{Shvarts:1984,Voronov:1984} that states that we only need probe things with Grassmann algebras, so supermanifolds of the form $\R^{0|p}$. Schwarz and Voronov refer to these `generalised points' as $\Lambda$-points.  In practice one can pick a `large enough' Grassmann algebra to perform the required calculations, but as nothing should really depend on this, one should check the functorial nature of all following constructions.
\begin{remark}
With no details, there is an amazing link between the De-Witt--Roger's approach to supermanifolds and the locally ringed space approach via $\Lambda$-points. Very loosely, the set of $\Lambda$-points for a given Grassmann algebra is a smooth manifold (with an extra structure on its tangent bundle), and smooth functions on it are smooth functions in the sense of Rogers.   This was probably first spelt out by Schmitt \cite{Schmitt:1997} who was inspired by the more categorical approach of  Molotkov \cite{Molotkov:1984}. I would not be surprised if this was know much earlier.  For details consult Schmitt.  
\end{remark}
In more generality, fields whose targets are supermanifolds (irrespective of the source) need to be thought of in terms of the internal homs as parametrised functors. We will see this again.

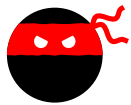
\begin{center}
\Ninja [5]
\end{center}

\section{Space-time and spinors}
\subsection{Recap} 
Via supermanifolds and generalised supermanifolds we have made  sense of anticommuting degrees of freedom, although this is in terms of functors and natural transformations. Anticommuting fields are essential to describe physical fermions (e.g., the electron) in a classical setting.  We will need to discuss spin and spinors if we really want to discuss fermions\par

\subsection{Minkowski space-time and the Poincar\'{e} group}
Maybe the most important group in physics is the \emph{Poincar\'{e} group} as the positive energy unitary irreducible representations are indexed by mass (non-negative number) and spin (integer or half integer), and classify  particles in quantum field theory. Thus, looking at representations of the Poincar\'{e} group, that is the group of isometries of Minkowski space-time, is essential in physics -- we need \emph{spinors} to understand the spin  $\half$-particles. We will present the most direct route to spinors understood as things that \emph{Clifford algebras} act on, and state that they do carry a representation of the Lorentz group (via the associated Lie algebra). This is only part of the much richer story. \par
Einstein took the idea of unifying space and time into one geometric entity very seriously  and this lead to special relativity and general relativity. We will only discuss special relativity in these lectures -- that is we will forget gravity completely. We will work with the `mostly positive' conventions for reasons that will be explained.\par 
\emph{Minkowski space-time} we will understand as the pseudo-Riemannian manifold $\R^{3,1} : = \big ( \R^{4}, \eta \big )$ where the metric in rectangular coordinates (`inertial coordinates') is 
$$\eta_{\mu \nu} = \textnormal{diag} (-1, + 1, + 1 + 1)\,.$$
The \emph{space-time interval} (infinitesimal) is defined as %

$$ - \delta s^2  = - c^2 \delta t^2 + \delta x^2 + \delta y^2 + \delta z^2\,,$$
where $c$ is the speed of light -- we will set $c = 1$ from now on. All inertial observers will agree on the space-time interval between events, even if they don't agree on the duration or distance. Everything is really encoded in the geometry of Minkowski space-time. In particular, the world line of light rays always satisfies $\delta s^2 =0 $ and this encodes the causal structure of special relativity. \par 
\begin{remark} We can work with Minkowski space-time of different spacial dimensions, for example in string theory $9+1$ dimensions is needed. We will  exclusivly stick to $3+1$ in these lectures and just remark that mod some technicalities with gamma matrices and spinors (see the next subsection)  the ideas remain the same in any dimension.  
\end{remark}
The \emph{Poincar\'{e} group} is the group of isometries of Minkowski space-time
$$O(3,1) \rtimes \R^4\,,$$
which consists of $6$ rotations + boost  and $4$ translations -- we have a $10$ dimensional noncompact Lie group. The \emph{Lorentz group} is then the stabiliser subgroup of the origin, i.e., forget translations.\par 
Infinitesimally, the Killing vector fields on Minkowski space-time of course form a Lie algebra -- the \emph{Poincar\'{e} algebra}. In rectangular coordinates the Killing equation is given by
$$\partial_{\mu}X_{\nu} + \partial_{\nu}X_{\mu} =0\,.$$
Taking the derivative again and cycling ober all indices, one can see that a Killing vector field must be of the form
$$X_{\mu} = x^{\nu}b_{\mu } + a_{\mu}\,.$$
Substituting this back into the Killing equation we observe that $b_{\mu \nu} = - b_{\nu \mu}$ and so we have $10$ independent vector fields given by choices of $b$ and $a$. A convenient choice is to first consider $b=0$, and then $a=0$ and $b_{\mu \nu} = 1$ when $\mu > \nu$. This leads to 
\begin{align}
& P_{\mu} = \frac{\partial}{\partial x^{\mu}},\\
& J^{\mu \nu} = x^{\mu}\eta^{\nu \lambda}\frac{\partial}{\partial x^{\lambda}} \: {-} \:   x^{\nu}\eta^{\mu \lambda}\frac{\partial}{\partial x^{\lambda}}.
\end{align}
One can then show that the Poincar\'{e} algebra can be written as
\begin{align}\label{eqn:LorentzAlgebra}
\nonumber & [P_{\mu}, P_{\nu}] =0,\\
\nonumber & [P_{\mu} ,J^{\lambda \sigma}] = \big (\delta^{\lambda}_{\mu} \eta^{\sigma \rho}  \: {-} \: \delta^{\sigma}_{\mu} \eta^{\lambda \rho}  \big)P_{\rho},\\
& [J^{\mu \nu}, J^{\rho \sigma}] = \eta^{\nu \rho}J^{\mu \sigma } \: {-} \: \eta^{\mu \rho}J^{\nu \sigma } \: {-}\: \eta^{\nu \sigma}J^{\mu \rho} \:{+}\: \eta^{\mu \sigma}J^{\nu \rho }.
\end{align}
Physically relevant fields will have nice transformations properties under the Poincar\'{e} group, or infinitesimally the fields will carry a representation of the  Poincar\'{e} algebra -- really just means we know how they transform under infinitesimal Lorentz transformations and translations. Amazingly, the (physically interesting) representations are not exhausted by tensor representations, we also have spinor representations.

\subsection{Clifford algebras}
Let $V$ be a real vector space equipped with a symmetric bilinear form $B(- , - )$. The \emph{Clifford algebra} associated with $(V,B)$, which we denote as $Cl(V)$ assuming no ambiguity in the bilinear form, is defined as follows. First, take the free tensor algebra
$$\bigotimes V : = \oplus_{j=0}^{\infty} V^{\otimes j}\,,$$
where we take $V^{\otimes 0} = \R$. Then consider the ideal $\mathcal{I}$  in $\bigotimes V$ spanned by
$$\{ x\otimes y + y \otimes x  + 2 B(x,y) \: | \: x,y \in V\} \,.$$
By picking a basis $(e_{1}, \cdots e_{2}) $ of $V$, then as an algebra $Cl(V)$ is generated by $(\Id, e_{1}, \cdots e_{k}) $ subject to the `Clifford condition' (we drop $\otimes$ from the notation)
$$e_i e_j + e_j e_i = - 2 B(e_i, e_j)\,.$$
\begin{example}
If $B=0$ then we have a Grassmann algebra. Thus in a loose sense, Clifford algebras are the quantisation of Grassmann algebras.
\end{example}
\begin{example}
Consider $V = \R$ and equipped with a normed basis $\epsilon$ (i.e., $\langle \epsilon , \epsilon \rangle = 1 $). Then $Cl(\R) \simeq \{a + \epsilon b \: | \: a,b \in \R\} $ as a vector space. However, as $\epsilon^2 = -1$ from the Clifford condition, as algebras $Cl(\R) \simeq \C$.
\end{example}

 \subsection{The Clifford algebra of space-time and spinors}
We can use the Minkowski metric on $\R^{3,1}$ to define a symmetric bilinear form and then employ an orthonormal basis to slightly simplify the Clifford condition.  In doing so we obtain the \emph{Clifford--Dirac relation}
 $$\{ \gamma^{\mu} , \gamma^{\nu} \} :=  \gamma^{\mu} \gamma^{\nu}  + \gamma^{\nu} \gamma^{\mu}  = 2 \eta^{\mu \nu} \Id\,.$$
The crucial observation here is that we can build an abstract representation of the  \emph{Lorentz algebra}.  If we define $\gamma^{\mu \nu} :=  \gamma^{[\mu} \gamma^{\nu ]} =  \half( \gamma^{\mu} \gamma^{\nu }  -  \gamma^{\nu} \gamma^{\mu })$, then  $\Sigma^{\mu \nu} = \half\gamma^{\mu \nu}$ satisfies (we won't give a proof here)
\begin{equation}\label{eqn:SpinRep}
[\Sigma^{\mu \nu}, \Sigma^{\rho \sigma}] = \eta^{\nu \rho}\Sigma^{\mu \sigma} -  \eta^{\mu \rho}\Sigma^{\nu \sigma} + \eta^{\mu \sigma}\Sigma^{\nu \rho} - \eta^{\nu \sigma}\Sigma^{\mu \rho}\, ,
\end{equation}
which you should then compare with (\ref{eqn:LorentzAlgebra}). This is vital in understanding the importance of spinors in physics.\par 
 This is of course rather abstract at this stage, one then usually chooses matrix  representations of the $\gamma$'s to build a matrix representation of the Clifford algebra, which we denote as $Cl(3,1)$. There are many  representations that one could choose -- we will pick \emph{real gamma matrices} as we are interested in geometry and would like to avoid complex numbers at this stage. This is also reflected in our conventions with regards to the signature of the Minkwoski metric: we cannot pick real gamma matrices  with the `mostly negative' conventions. As such we  will choose the following \emph{real Majorana representation} 
\[ \gamma^{0} = \left( \begin{array}{rrrr}
0 & +1 & 0 & 0 \\
-1 & 0 & 0 & 0  \\
0 & 0 & 0  & -1\\
0 & 0 & +1  & 0\\
\end{array} \right), \hspace{25pt}
 \gamma^{1} = \left( \begin{array}{rrrr}
0 & +1 & 0 & 0 \\
+1 & 0 & 0 & 0  \\
0 & 0 & 0  & +1\\
0 & 0 & +1  & 0\\
\end{array} \right),\]
\\
\[ \gamma^{2} = \left( \begin{array}{rrrr}
+1 & 0 & 0 & 0 \\
0 & -1 & 0 & 0  \\
0 & 0 & +1  & 0\\
0 & 0 & 0  & -1\\
\end{array} \right), \hspace{25pt}
 \gamma^{3} = \left( \begin{array}{rrrr}
0 & 0 & 0 & +1 \\
0 & 0 & -1 & 0  \\
0 & -1 & 0  & 0\\
+1 & 0 & 0  & 0\\
\end{array} \right).\]
These $4\times 4$ real matrices act on $\R^{4}$ in the usual way. With this in mind, we adopt the convention of spin indices of $(\gamma^{\mu})_{\alpha}^{\:\:\: \beta}$. We then define \emph{Majorana spinors} as the ``things" that these gamma matrices act on.  We thus write $\mathbb{M} = \R^{4}$ for the space of all Majorana spinors. We will conventionally choose the components of a Majorana spinor to have lower spin indices, and so we write $u_{\alpha}$.  \par 
As Majorana spinors are ``things'' that $Cl(3,1)$ act on, then because of (\ref{eqn:SpinRep}) we have an action of the Lorentz Lie algebra on spinors -- we think of this as an action of infinitesimal Lorentz transformations. The specification of any infinitesimal Lorentz transformation (3 infinitesimal rotations and 3 infinitesimal boosts) is encoded in an antisymmetric tensor $\omega_{\mu \nu} = - \omega_{\nu \mu}$. One can then write
$$\delta_{\omega} u_{\alpha}  = \frac{1}{4} \omega_{\mu \nu} (\gamma^{\mu \nu})_{\alpha}^{\:\: \: \beta}u_{\beta}\, ,$$ 
where the numerical factor comes from a factor of $\half$ due to the antisymmetric properties of the tensors involves -- to avoid over-counting.\par 
By taking the exponential we obtain a finite Lorentz transformation (as we integrate a Lie algebra we only get the component connected to the identity, though this need not concern us)
$$u'_{\alpha} = \exp( \frac{1}{4} \omega_{\mu \nu}\gamma^{\mu \nu} )_{\alpha}^{\:\:\: \beta} u_{\beta}\,.$$
\begin{remark}
\emph{Dirac spinors} are understood as the complexification of the Majorana spinors, $\mathbb{D} :=  \mathbb{M} \otimes_{\R} \C \simeq \C^4$. \emph{Weyl spinors} are then defined in terms of a decomposition of the Dirac spinors into a direct sum $\mathbb{D} = \mathbb{W}_{+}\oplus \mathbb{W}_{-}$ defined by chirality.  More explicitly, we first define $\gamma_{5}:= \gamma_{0}\gamma_{1}\gamma_{2}\gamma_{3}$ (where we have used the Minkowski metric to lower the indices). The decomposition is then defined in terms of the eigenvalues of $\gamma_{5}$. That is,  $u \in \mathbb{W}_{\pm}$ if  $\gamma_{5} u = \pm \rmi \: u$. Via this decomposition we see that $\mathbb{W}_{\pm} \simeq \C^2$.
\end{remark}
As defined here, the Majorana spinor as \emph{commuting objects}, they are  described by the coordinates on $\mathbb{M}$ considered as a linear manifold. In order to have \emph{anticommuting Majorana spinors} -- which are more common in physics -- we need to employ the parity reversion functor. Thus, we consider anticommuting Majorana spinors to be described by the coordinates on $\Pi \mathbb{M} \simeq \R^{0 | 4}$.\par 
The manifold $\mathbb{M}$ (and/or the supermanifold $\Pi \mathbb{M}$) comes equipped with  the \emph{charge conjugation tensor} --  this operator exchanges particles and antiparticles. The defining property is
$$C \gamma^{\mu}C^{-1} = - (\gamma^{\mu})^{\textnormal{t}}\,.$$
Or written more explicitly
$$  C^{\alpha \gamma}(\gamma^{\mu})_{\gamma}^{\:\: \delta} C_{\delta \beta} = {-}(\gamma^{\mu})^{\alpha}_{\:\: \beta}
\,,$$
Thus, considered as a matrix,  in our chosen representation, $ C =  -\gamma^{0}$. To set some useful notation for later, we define $(C\gamma^{\mu})^{\alpha \beta} :=  C^{\alpha \delta }(\gamma^{\mu})_{\delta}^{\:\: \beta}$, and a direct computation show that the resulting tensor is symmetric.  
\subsection{Anticommuting Majorana spinor fields}
All fields in physics are understood as sections of the appropriate fibre bundles over space-time -- in the setting of special relativity we can usually take these fibre bundles as being trivial. Anticommuting Majorana spinor fields we understand as sections of the trivial (vector) bundle whose typical fibre is $\Pi \mathbb{M}$, but we need to modify our notion of sections to capture the anticommuting nature of these fields.\par 
Consider the trivial bundle $\pi :  \R^{3,1} \times \Pi \mathbb{M} \longrightarrow \R^{3,1}$. For notational simplicity let us set $E :=  \R^{3,1} \times \Pi \mathbb{M}$. Then consider the pullback bundle 

\begin{center}
\leavevmode
\begin{xy}
(0,20)*+{\pi_{M}^{*}E}="a"; (30,20)*+{E}="b";%
(0,0)*+{S \times \R^{3,1}}="c"; (30,0)*+{\R^{3,1}}="d";%
{\ar "a";"b"}?*!/_3mm/{\hat{\pi}};
{\ar "a";"c"}?*!/^7mm/{(\Id_{S}, \pi)};
{\ar "b";"d"}?*!/_5mm/{\pi};
{\ar "c";"d"}?*!/^3mm/{\pi_{M}};
\end{xy}
\end{center}
 We have used the canonical identification $\pi_{M}^{*}E =  S \times E$ in the above diagram. We then define Majorana spinor fields (from now we always assume spinor fields are anticommuting unless we explicitly say otherwise) parameterised by $S \in \catname{SM}$ to be sections of the above pullback bundle. More carefully, the generalised supermanifold of (global) Majorana spinor fields is the functor 
$$\underline{\Sec}(E)(-) : \catname{SM}^\opp \longrightarrow \catname{Set}$$
defined as
$$\underline{\Sec}(E)(S)  :=  \left \{  \Psi_{S} \in \Hom(S \times \R^{3,1} , E) \:  | \: (\Id_{S}, \pi )\circ \Psi_{S}  = \pi_{M} \right \},$$
and for any $\phi \in \Hom(P,S)$ 
\begin{align*}
\underline{\Sec}(E)(\phi) & : \underline{\Sec}(E)(S) \longrightarrow \underline{\Sec}(E)(P)\\
& \Psi_{S} \mapsto \Psi'_{P} := \Psi_{S} \circ (\phi , \Id_{M}). 
\end{align*}
A little more explicitly, we can employ global coordinates  $(x^{\mu}, \psi_{\alpha})$  on $R^{3,1} \times \Pi \mathbb{M}$,  and using  the theorem of Schwarz and Voronov \cite{Shvarts:1984,Voronov:1984} we can set $S = \Lambda :=\R^{0|p}$ without loss of generality. Let us then employ global coordinates $\zx^{i}$ on $\R^{0|p}$, and then we can write
$$\big (\zx^i , x^\mu , \psi_{\alpha} \big) \circ \Psi_{\Lambda} = \big (\zx^i , x^\mu , \psi_{\alpha}(\zx, x) \big ).$$
It is customary to simply write
$\psi_{\alpha}(\zx, x) := \psi_{\Lambda \alpha}(x)$, or drop reference to the Grassmann algebras (or more general supermanifolds) altogether. \par 
The upshot is that we can not really think of a single Majorana spinor field, but rather we need to think in terms of a parametrised family of maps -- such fields are best thought of as functors from $\catname{SM}^\opp$ to $\catname{Set}$ that are parametrised by space-time.  
  
\subsection{An Aside: The Dirac and Klein--Gordon equations}
Spinors in physics were discovered by Dirac (1928) via `taking the square root' of the \emph{Klein--Gordon equation} (1926). Later, it was discovered that mathematicans, and of course Clifford himself,  had already considered such things. To sketch, the relation let us start from the \emph{Dirac equation}
$$(\gamma^{\mu} \partial_{\mu} + m)_{\alpha}^{\:\: \beta} \psi_{\beta}(x) =0\,,$$
we have already sketched what $\psi_{\beta}(x)$ means, really we should be thinking of a family of Dirac equations parameteried by (at least) all finite dimensional Grassmann algebras. Then applying the operator on the left again, but this time with negative mass (the sign of the mass term is irrelevant in the Dirac equation) we obtain
$$(\gamma^{\mu} \partial_{\mu} - m)_{\alpha}^{\:\: \beta}   (\gamma^{\nu} \partial_{\nu} + m)_{\beta}^{\:\: \gamma} \psi_{\gamma}(x) =0\,.$$
Using the Clifford--Dirac relation we obtain
$$(\eta^{\mu \nu} \partial_{\nu} \partial_{\mu} - m^2)\psi_{\alpha}(x)=0\,,$$
or in more traditional notation
$$\left({-}\frac{\partial^2}{\partial t^2} +  \nabla^2\right)\psi_{\alpha}(x) = m^2 \psi_{\alpha}(x)\,.$$
That is, if a spinor satisfies the Dirac equation then each of the components satisfies the Klein--Gordon equation.
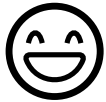
\begin{center}
\Laughey [5]
\end{center}
\section{Real super-Minkowski space-times}
\subsection{Recap}
We have now gone from some mechanical motivations, through supermanifolds and then Clifford algebras and spinors to arrive at a position from which we can look at a version of super-Minkowski space-time -- the geometric arena of supersymmetry. At the fundamental level, one  extends Minkowski space-time by appending anticommuting spinor coordinates in order to geometrically realise the supersymmetry algebra. 

\subsection{The $N=1$ super Poincar\'{e} algebra}
We extend the  Poincar\'{e}  algebra by adjoining anticommuting  \emph{Majorana spinor} generators $Q^{\alpha}$ (see \cite{Golfan:1971}).  The resulting algebra under study is thus
\begin{subequations}
\begin{align}
& [P_{\mu} , J^{\lambda \sigma}] = (\delta_{\mu}^{\:\: \lambda} \eta^{\sigma \rho} -\delta_{\mu}^{\:\: \sigma} \eta^{\lambda \rho} )P_{\rho},\\
&[J^{\mu \nu} , J^{\rho  \sigma}] = \eta^{\nu \rho}J^{\mu \sigma} - \eta^{\mu \rho}J^{\nu \sigma}  - \eta^{\nu \sigma}J^{\mu \rho} + \eta^{\mu \sigma}J^{\nu \rho},\\
& [J^{\mu \nu}, Q^{\alpha}]  = - \frac{1}{4}Q^{\beta}(\gamma^{\mu \nu})_{\alpha}^{\:\: \beta},\\
& [Q^{\alpha}, Q^{\beta}] = \frac{1}{2}  (C\gamma^{\mu})^{\alpha \beta}P_{\mu} \,. \label{Eqn:4DSUSY}
\end{align} 
\end{subequations}
All other Lie brackets are zero. The resulting algebra contains the Poincar\'{e}  algebra as the degree zero/`commuting' component.  Note that we have considered everything to be $\Z_2$-graded and thus the bracket of two $Q$'s is an anticommutator in `traditional' language.\par 
The first two equations are simply the standard  Poincar\'{e}  algebra. In particular, the second equation is the Lorentz subalgebra. The third equation tells us that the generators are spinors -- this defines how they transform under infinitesimal Lorentz transformations. The third equation is the important one from our perspective -- it show that the supersymmetry generators $Q$ are the `square root' of the 4-momentum operator, i.e., the generator of translations in Minkowski space-time. It is this fourth equation \eqref{Eqn:4DSUSY} that we we focus on.

\subsection{Constructing super-Minkowski space-time}
The idea is to find a supermanifold that allows us to geometrically as vector fields realise the $N=1$ super Poincar\'{e} algebra -- or really we will focus on just realising $Q$'s and $P$, that is we will forget about the Lorentz algebra. One can make educated guesses into what the supermanifold should be and then the form of the vector fields. We will however, take a slightly more formal route.\par 
Recall that given some Lie algebra we can integrate this to obtain a Lie group. The ``Lie group--Lie algebra correspondence'' can be made quite explicit using the Campbell--Baker--Hausdorff formula, and this is especially true when we have a \emph{nilpotent Lie algebra}, i.e.,  the nested Lie brackets varnish at some stage. In particular, there are no questions about convergence. The construction gives a local form (near the identity) of the group multiplication in coordinates. \par 
We claim -- without proof -- that the same method can be applied to the case at hand. To make real sense of this one should appeal to the functor of points in order to make this rigorous, however we will simply proceed formally.\par 
With these comments in mind, let us consider two general elements in the super translation algebra (forgetting the Lorentz part)
\begin{align*}
A =  x^{\mu} P_{\mu} + \theta_{\alpha} Q^{\alpha}, & & B =  x'^{\mu} P_{\mu} + {\theta'}_{\alpha} Q^{\alpha}\,,
\end{align*}
Note that $x$ are taken to be even, while $\theta$ are taken to be odd. Thus, a general element we see is an even object. Then using the Campbell--Baker--Hausdorff formula, noting that we only need to consider terms up to one Lie bracket,  we obtain
\begin{equation}\label{eqn:GroupLaw}
 A\circ B  = A+B + \half [A,B] =  \left( x^{\mu} + x'^{\mu} + \frac{1}{2} \theta'_{\beta} \theta_{\alpha}(C\gamma^{\mu})^{\alpha \beta}  \right )P_{\mu} + (\theta_{\alpha} + \theta'_{\alpha})Q^{\alpha}\,.
 \end{equation}
From this one can formally read off the group structure in terms of coordinates, which justifies our very suggestive notation. The upshot of this analysis is that we have uncovered the basic structure of super-Minkowski space-time.\par 
\begin{definition}
\emph{Super-Minkowski space-time}  $M^{3,1 | 4}$ is the supermanifold that can be equipped with global coordinates
$$(x^{\mu} , \theta_{\alpha})\,,$$
where the anticommuting coordinates are Majorana spinors. 
\end{definition}
Thus, given how we have defined  Majorana spinors we have the identification (via a choice of coordinates)
$$M^{3,1 | 4} \simeq \R^{4|4} \simeq \R^{4}\times \Pi \mathbb{M} \,.$$
Via construction, the reduced manifold is given by standard Minkowski space-time.\par 
From the (formal) group law given by (\ref{eqn:GroupLaw}) the supersymmetry transformations on $M^{3,1 | 4}$ are given by 
\begin{subequations}
\begin{align}
 & x^{\mu} \mapsto x^{\mu} + \frac{1}{4}\epsilon_{\beta}\theta_{\alpha} (C\gamma^{\mu})^{\alpha \beta}\,,\\
 & \theta_{\alpha} \mapsto \theta_{\alpha} + \epsilon_{\alpha}\,,
\end{align}
\end{subequations}
where $\epsilon_{\alpha}$ is a Majorana spinor-valued anticommuting parameter. 

\subsection{The SUSY generators and covariant derivatives}
The goal is now to represent $Q$ and $P$ as vector fields on  super-Minkiowski space-time. From the group law we can proceed formally and write down the left invariant vector fields -- this is just given by taking the derivative with respect to the primed variables. In this way we obtain
\begin{align*}
& P_{\mu} = \frac{\partial}{\partial x^{\mu}}, && Q^{\alpha} = \frac{\partial }{\partial \theta_{\alpha}} + \frac{1}{4}\theta_{\beta}(C\gamma^{\mu})^{\beta \alpha}\frac{\partial}{\partial x^{\mu}}\,.
\end{align*} 
By `hand' you can easily check that these satisfy the correct commutation rules.\par 
 As in the mechanical case, we need to construct  covariant derivatives if we are really to work with supersymmetry. The method is quite simple, we look at the right invariant vector fields -- this then comes down to taking the derivative of the group law with respect to the unprimed variables. In this way we obtain (and it might have been your first guess)
 $$\mathbb{D}^{\alpha} =  \frac{\partial }{\partial \theta_{\alpha}} {-} \frac{1}{4}\theta_{\beta}(C\gamma^{\mu})^{\beta \alpha}\frac{\partial}{ \partial x^{\mu}}\,.$$

\subsection{The SUSY structure}
From a geometric point of view the distribution spanned by the SUSY covariant derivatives is an object of central interest -- as far as I know it was Manin who first took this point of view (see \cite{Manin:1991}).  Interestingly, this distribution is \emph{maximally non-integrable}. \par 
Let us quickly recall the notion of maximally non-integrable distributions. A distribution is a subsheaf of the tangent sheaf that is locally a direct factor -- we need not elaborate on this too carefully as we will have global structures. Let us decompose the module of vector fields on super-Minkowski space-time as a direct sum of a distribution that we want to study, and another module of whatever is not in this distribution
$$\Vect(M^{3,1 | 4}) = \mathcal{D} \oplus \mathcal{N}\,,$$
with $\mathcal{N} := \Vect(M^{3,1 | 4})\setminus \mathcal{D}$. We define map
$$\phi : \Vect(M^{3,1 | 4})  \longrightarrow \Vect(M^{3,1 | 4})/ \mathcal{D} =: \mathcal{N}\,, 
$$
which is `remove anything in the distribution'. The \emph{Frobenius curvature} is then given by
\begin{align*}
R &: \mathcal{D} \times \mathcal{D} \longrightarrow \mathcal{N}\\
 & (X,Y) \mapsto \phi([X,Y])\,.
\end{align*}
Basically, one is looking how much Lie brackets of elements in the distribution lie outside the distribution. Of course, if the distribution is closed with respect to the Lie bracket then we say that the distribution is \emph{integrable}. The characteristic distribution $\mathcal{C}$ is defined as all $X \in \mathcal{D}$ such that $R(X, -) =0$. If $\mathcal{C}$ consists of only the zero vector field then $\mathcal{D}$ is \emph{maximally non-integrable}.  This means that all Lie brackets of non-zero elements of $\mathcal{D}$ fall outside of $\mathcal{D}$. Nothing in the above is particular to supergeometry or supersymmetry. \par 
Now, let us specialise to the case of 
$$\mathcal{D} = \Span \left \{\mathbb{D}^{\alpha}  \right\}\,.$$
Clearly, this a distribution -- you can check that any vector field on super-Minkwoski space-time can be written using $P_{\mu}$ and $\mathbb{D}^{\alpha}$ as a basis. Moreover, the rank of $\mathcal{D}$ is $(0|4)$, i.e., has a basis consisting of four odd vector fields. In slightly more generality, Manin calls such distributions \emph{pre-SUSY structures} -- we have an odd basis that is of maximal dimension. If this distribution is maximally non-integrable then we have a \emph{SUSY structure}.\par
Directly  we examine the Frobenius curvature by looking at $X = X_{\alpha}(x, \theta)\mathbb{D}^{\alpha}$ and $\mathbb{D}^{\beta}$, this is enough to get at the characteristic distribution.  Thus,
\begin{align*}
  R([X, \mathbb{D}^{\beta}] ) & = \phi \left(X_{\alpha}[\mathbb{D}^{\alpha}, \mathbb{D}^{\beta}] \pm \big(\mathbb{D}^{\beta}X_\alpha\big)\mathbb{D}^\alpha \right)\\ 
  & = \frac{1}{2}X_{\alpha} (C\gamma^{\mu})^{\alpha \beta} \frac{\partial}{\partial x^{\mu}}\,.
\end{align*}
Then we see that this is zero if and only if $X_{\alpha} =0$, thus the characteristic distribution consists of only the zero vector and so $\mathcal{D}$ is maximally non-integrable.
\begin{remark}
In the one dimensional case, i.e., mechanics, the distribution spanned by the SUSY covariant derivatives is still maximally non-integrable and this time of corank $(1|0)$. Thus we have a (super generalisation of a)  \emph{contact structure}. SUSY structures are a particular higher dimensional generalisation of contact structures. 
\end{remark}

\subsection{Aside: Chiral Superspace}
In practice, one often uses Weyl spinors instead of Majorana spinors (or their complexification, the Dirac spinors) when building supersymmetric field theories. Without details, the basic idea remains the same -- one extends space-time by appending some anticommuting spinor coordinates, examines the associated supersymmetry transformations and build the generators and covariant derivatives. Modifications of the methods presented in the Section \ref{Sec:mechanics} allow one to build supersymmetric actions from superfields and covariant derivatives (often one needs further constraints on the superfields). All this is outside the scope of this first look at supersymmetry.
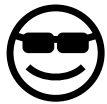
\begin{center}
\Cooley[5]
\end{center} 
 %
\section{Round-up and Conclusion}
In these lectures we have motivated some of the mathematics behind supersymmetry and in particular the use of supermanifolds. One aspect that we have explained is the nature of Grassmann odd fields as parametrised functors. This point of view is not usually stressed in the physics literature. From there we build odd spinor fields and defined super-Minkowski space-time as well as geometrically defining supersymmetry. Introductory texts on supersymmetry often stress the operational use of superspace but do not properly explain the language pertaining to differential geometry. It became clear to me during my Master's degree, while I was totally unaware of supergeometry, that supersymmetry should be formulated in a geometric way.  The categorical approach to supermanifolds and odd maps is barely touched upon in the physics literature. Everything is understood formally or operationally and the attitude is that physicists need not worry about the exact nature of the objects pushed around by pencil on paper.  \par 
We have not attempted to build interesting and phenomenologically realistic models using superspace methods. Nor have we examined any quantum models. The aims of the lectures have been modest and provide a geometric glance at supersymmetry, albeit one that forces us to confront ideas in category theory early on.

\subsection*{Further reading}
\begin{description}
\item[C. Carmeli, L. Caston \& R. Fioresi]  \href{http://www.ems-ph.org/books/book.php?proj_nr=132}{Mathematical foundations of supersymmetry}, EMS Series of Lectures in Mathematics, \emph{European Mathematical Society} (EMS), Z\"{u}rich, 2011, xiv+287 pp. ISBN: 978-3-03719-097-5.

\item [D.A. Le\u{\i}tes]
Introduction to the theory of supermanifolds,
 \emph{Uspekhi Mat. Nauk}, \textbf{35}(1(211)):3--57, 255, 1980.
 
\item  [Y.I. Manin]
{ Gauge field theory and complex geometry}, volume 289 of {\em
  Grundlehren der Mathematischen Wissenschaften [Fundamental Principles of
  Mathematical Sciences]}.
 Springer-Verlag, Berlin, second edition, 1997.
Translated from the 1984 Russian original by N. Koblitz and J. R.
  King, With an appendix by Sergei Merkulov.
  
\item [V.S. Varadarajan] \href{http://bookstore.ams.org/cln-11}{Supersymmetry for mathematicians: an introduction},
Courant Lecture Notes in Mathematics, 11. New York University, Courant Institute of Mathematical Sciences, New York; \emph{American Mathematical Society}, Providence, RI, 2004. viii+300 pp. ISBN: 0-8218-3574-2.

\item[S. Weinberg] \href{https://www.cambridge.org/ae/universitypress/subjects/physics/theoretical-physics-and-mathematical-physics/quantum-theory-fields-volume-3}{The quantum theory of fields. Vol. 3: Supersymmetry}, \emph{Cambridge University Press}, 2000, ISBN:9781139632638.

\item[J. Wess \& J. Bagger]
 \href{https://press.princeton.edu/books/paperback/9780691025308/supersymmetry-and-supergravity}{Supersymmetry and Supergravity}, \emph{Princeton University Press}, Princeton, 1992, ISBN 0-691-02530-4.
 
 \item[P. West]
 \href{9789813103719, 981310371X}{Introduction To Supersymmetry And Supergravity} (Revised And Extended 2nd Edition), \emph{Singapore: World Scientific Publishing}, 1990, ISBN 978-981-02-0098-5.
\end{description}

\appendix
\section{Appendix}
\subsection{Lie Superalgebras}\label{appendix:Lie}
A Lie superalgebra is a $\Z_2$-graded generalisation of a Lie algebra. Suppose we have a $\Z_2$-graded vector space (over $\R$ or $\C$, say) 
$$\mathfrak{g} = \mathfrak{g}_0 \oplus \mathfrak{g}_1\,,$$
and that we define the parity of elements as $\widetilde{A} \in \Z_2$ depending on what subspace the element $A$ belongs to. In general, an element is the sum of both an even and odd element, and we extend the following to inhomogeneous elements via linearity.  
A \emph{Lie superalgabra} is a  $\Z_2$-graded vector space equipped with a bilinear product
$$[-,-] : \mathfrak{g} \times \mathfrak{g} \rightarrow \mathfrak{g}\,,$$
 called the \emph{Lie bracket} that satisfies the following
 \begin{enumerate}
 \item $\widetilde{[A,B]} = \widetilde{A}+ \widetilde{B}$,
 \item $[A,B] = - (-1)^{\widetilde{A} \widetilde{B}}\, [B,A]$,
 \item $[A,[B,C]] = [[A,B],C]+ (-1)^{\widetilde{A} \widetilde{B}}\, [B,[A,C]]$.
 \end{enumerate}

 The third property is called the \emph{Jacobi identity}.

\subsection{Grassmann Algebras}
\begin{definition}
 A \emph{Grassmann algebra} (over $\R$) with $q$ generators  $\Lambda(\zx^1, \cdots ,\zx^q)$ consists, of a polynomials generated by the  $\{ \zx^\alpha \}$, subject to the  relation
$$\zx^\alpha \zx^\beta = -  \zx^\beta \zx^\alpha\,.$$
\end{definition}
Note that the above relation implies that $(\zx)^2 =0$ and thus the maximal order of a polynomial in  $\Lambda(\zx^1, \cdots ,\zx^q)$ is $q$. Specifically, a given generator can only appear at most once in a given monomial. Grassmann algebras are $\Z_2$-graded, i.e.,
$$\Lambda(\zx^1, \cdots ,\zx^q) = \Lambda_0(\zx^1, \cdots ,\zx^q)\oplus\Lambda_1(\zx^1, \cdots ,\zx^q)\,.$$

The even elements $f\in \Lambda_0(\zx^1, \cdots ,\zx^q)$ are polynomials consisting of monomials of even order, and odd elements $g\in \Lambda_1(\zx^1, \cdots ,\zx^q)$ are polynomials consisting of monomials of odd order. Let us define the parity of elements as $\widetilde{f} \in \Z_2$ depending on what subspace the element $f$ belongs to. In general, an element is the sum of both an even and odd element, and we extend the following to inhomogeneous elements via linearity.\par 
Grassmann algebras are supercommutative superalgebras as
$$f\cdot g = (-1)^{\widetilde{f} \widetilde{g}}\,g \cdot f\,. $$

\subsection{Calculus with Anticommuting Variables}\label{App:calculus}
Differentiation with respect to Grassmann generator is defined algebraically. In particular, 
$$\frac{\partial \zx^\alpha}{\partial \zx^\beta} = \delta_\beta^\alpha\,.$$
As elements in a Grassmann algebra are polynomials, the above definition allows us to understand derivatives of general elements of the algebra. One has to be careful with minus signs.  For example, as $\zx^\alpha \zx^\beta = -  \zx^\beta \zx^\alpha$, we observe that 
$$\frac{\partial \zx^\alpha \zx^\beta}{\partial \zx^\gamma} = \delta_\gamma^\alpha \zx^\beta - \zx^\alpha \delta_\gamma^\beta\,. $$ 
Similar reasoning shows that 
$$\frac{\partial}{\partial \zx^\alpha} \frac{\partial}{\partial \zx^\beta} = -  \frac{\partial}{\partial \zx^\beta} \frac{\partial}{\partial \zx^\alpha}\,. $$
Note that this implies that $\big( \frac{\partial}{\partial \zx}\big)^2 =0$ for any given generator.\par 
The integration of odd variables is similarly defined algebraically in the sense that  is not a Lebesgue integral.
\begin{definition}
The \emph{Berezin integral} is the unique linear functional $\int D[\zx]-$ defined by the following properties
\begin{enumerate}
\item $\int D[\zx]\, \zx^1 \zx^2 \cdots \zx^q = 1$, and
\item $\int D[\zx]\, \frac{\partial f}{\partial \zx^\alpha} =0$,
\end{enumerate}
for any $f \in \Lambda(\zx^1, \cdots ,\zx^q)$.
\end{definition}
\begin{remark}
Other conventions with the normalisation appear in the literature.  Also note that there is a choice in how we order the generator of the Grassmann algebra and this can change overall signs.  
\end{remark}
Writing (formally) $D[\zx] = D[\zx^q]D[\zx^{q-1}] \cdots D[\zx^1]$, the Fubini law is
$$\int D[\zx]\,f(\zx) =  \int D[ \zx^q]\left(\cdots \left( \int D[\zx^2]\left(\int D[\zx^1]f(\zx^1,\zx^2, \cdots , \zx^q) \right) \right) \right)\,.$$
Any monomial that does not depend on some chosen $\zx^\alpha$ can be written as $w(\zx)= \frac{\partial}{\partial  \zx^\alpha}\big(\zx^\alpha \, w(\zx)\big)$. Thus, only highest order monomials are not automatically zero under integration.  Then for a general element $f(\zx) =  f_0 + \zx^1 f_1 + \cdots + \zx^1 \zx^2 \cdots \zx^q f_{q \cdots 21}$, the Berezin integral is
\begin{align*}
\int D[\zx]\,f(\zx)& = \int D[ \zx^q]\left(\cdots \left( \int D[\zx^2]\left(\int D[\zx^1]\zx^1 \zx^2 \cdots \zx^q f_{q \cdots 21} \right) \right)  \right)\\ &= \int D[\zx]\,\zx^1 \zx^2 \cdots \zx^q f_{q \cdots 21} = f_{q \cdots 21}\,.
\end{align*}
In short, upto normalisation and sign factors, the Berezin integral selects the coefficient of the top degree monomial of a general element of the Grassmann algebra. This implies the interesting relation that integration and differentiation are the same thing, i.e.,
$$\int D[\zx] =  \frac{\partial}{\partial \zx^q}\cdots \frac{\partial}{\partial \zx^2} \frac{\partial}{\partial \zx^1}\,.$$

\subsection{The Berezinian}\label{appendix:Ber}
Let us consider the Berezin integral on $\R^{p|q}$ defined in canonical coordinates $(x, \zx)$ (we suppress indices for notation ease),  $\int D[x,\zx]f(x, \zx)$. To evaluate this one first preforms the integration over the odd variable and then preforms the integration over the even variables as standard (here understood as the indefinite integral).  The question is what happens to the integral if we change variables, say $y = y(x, \zx)$ and $\theta = \theta(x, \zx)$ using the common abuse of notation. \par 
We define the \emph{Jacobian matrix} (with no justification here) as 
\[ J = 
\left(\begin{matrix}
   \frac{\partial y}{\partial x}       & - \frac{\partial y}{\partial \zx} \\
   \frac{\partial \theta}{\partial x}       &  \frac{\partial \theta}{\partial \zx} 
    \end{matrix}\right) = \left(\begin{matrix}
   A       & B \\
   C     &  D 
    \end{matrix}\right)\,.
\]  
We then define the Jacobian as the Berezinan (superdeterminant)of the Jacobian matrix, which is defined as
$$\Ber(J) := \det(A - BD^{-1}C)\det(D)^{-1}\,,$$
which is only defined if both block matrices $A$ and $D$ are invertible. The integration measures in the different coordinates are related by
$$ D[y, \theta] = D[x,\zx] \Ber(J)\,.$$
As a specific example, consider integration on $\R^{1|1}$ and the supersymmety transformations $t' =  t + \frac{1}{4} \epsilon  \theta$ and $\theta' = \theta + \epsilon$. Then 
\[ J = 
\left(\begin{matrix}
   \frac{\partial t'}{\partial t}       & - \frac{\partial t'}{\partial \theta} \\
   \frac{\partial \theta'}{\partial t}       &  \frac{\partial \theta'}{\partial \theta} 
    \end{matrix}\right) = \left(\begin{matrix}
   1      &\frac{1}{4}\epsilon \\
   0     &  1 
    \end{matrix}\right)\,,
\]  
and so $\Ber(J) = 1$. The integration measure is invariant under supersymmetry. \par 
Similarity, on $\R^{1,2}$ and the supersymmetry transformations  $t' = t +\frac{1}{4} \epsilon^{I} \theta^{J}\delta_{JI}$ and $ \theta'^{I} = \theta^{I} + \epsilon^{I}$, the Jacobian matrix can quickly be deduced to be
\[ J = 
 \left(\begin{matrix}
   1      &\star \\
   0     &  \Id_{2\times 2} 
    \end{matrix}\right)\,.
\]  

We do not actually care what the top right matrix is as the bottom left is the zero matrix. Clearly, $\Ber(J) =  1$, and again the integration measure in invariant under supersymmetry.


\begin{thebibliography}{10}


\begin{small}

\bibitem{Berezin:1975}
F. A. Berezin \& D.A. Le\u{\i}tes, Supermanifolds, (Russian)\emph{Dokl. Akad. Nauk SSSR} \textbf{224} (1975), no. 3, 505--508



\bibitem{Carmeli:2011}
C. Carmeli, L. Caston \& R. Fioresi,  \href{http://www.ems-ph.org/books/book.php?proj_nr=132}{Mathematical foundations of supersymmetry}, EMS Series of Lectures in Mathematics, \emph{European Mathematical Society} (EMS), Z\"{u}rich, 2011, xiv+287 pp. ISBN: 978-3-03719-097-5.

\bibitem{Castellani:2014}
L. Castellani, R. Catenacci \& P.A. Grassi,
Supergravity Actions with Integral Forms, \href{https://doi.org/10.1016/j.nuclphysb.2014.10.023}{\emph{Nucl. Phys. B}}, \textbf{889} (2014), 419--442. 

\bibitem{Castellani:2015}
L. Castellani, R. Catenacci \& P.A. Grassi,
The geometry of supermanifolds and new supersymmetric actions, \href{https://doi.org/10.1016/j.nuclphysb.2015.07.028}{\emph{Nucl. Phys. B}}, \textbf{899} (2015), 112--148. 

\bibitem{Castellani:2015b}
L. Castellani, R. Catenacci \& P.A. Grassi,
Hodge Dualities on Supermanifolds, \href{https://doi.org/10.1016/j.nuclphysb.2015.08.002}{\emph{Nucl. Phys. B}}, \textbf{899} (2015),  570--593.


\bibitem{Deligne:1999}
P.~Deligne \& D.S. Freed, 
Supersolutions, in \href{https://bookstore.ams.org/qft-1-2-s/}{\emph{Quantum fields and strings: a course for mathematicians}}, Vol. 1, 2 (Princeton, NJ, 1996/1997), 227--355,\emph{ Amer. Math. Soc., Providence, RI}, 1999. 

\bibitem{Gervais:1971}
J.L. Gervais \& B. Sakita,  Field theory interpretation of supergauges in dual models, \emph{Nuclear Physics B.} \textbf{34}(2) (1971), 632--639.

\bibitem{Golfan:1971}
Yu.A. Golfand \& E.P. Likhtman
Extension of the Algebra of Poincare Group Generators and Violation of p Invariance, \emph{JETP Lett.} \textbf{13} (1971) 323--326.

\bibitem{Green:1953}
H.S. Green,
A generalized method of field quantization, \emph{Phys. Rev.} \textbf{90} (1953), 270--273.

\bibitem{Freed:1999}
D.S.~Freed,
\href{https://bookstore.ams.org/fls/}{Five lectures on supersymmetry}, \emph{American Mathematical Society, Providence, RI}, 1999. viii+119 pp. ISBN: 0-8218-1953-4.

\bibitem{Leites:1980}
D.A. Le\u{\i}tes,
Introduction to the theory of supermanifolds,
 \emph{Uspekhi Mat. Nauk}, \textbf{35}(1(211)):3--57, 255, 1980.

\bibitem{MacLane1998}
S. Mac~Lane,
\newblock {\em Categories for the working mathematician}, volume~5 of {\em
  Graduate Texts in Mathematics}.
\newblock Springer-Verlag, New York, second edition, 1998.

\bibitem{Manin:1991}
Y.I. Manin,
{\em Topics in noncommutative geometry},
M. B. Porter Lectures, {\em Princeton University Press} (1991) 164 p.


\bibitem{Manin:1997}
Y.I. Manin,
{\em Gauge field theory and complex geometry}, volume 289 of {\em
  Grundlehren der Mathematischen Wissenschaften [Fundamental Principles of
  Mathematical Sciences]}.
Springer-Verlag, Berlin, second edition, 1997.

\bibitem{Molotkov:1984}
V. Molotkov, Banach supermanifolds, in Differential Geometric Methods in Theoretical Physics (Shumen,
1984), \emph{World Sci. Publishing}, Singapore, 117--125 1986.

\bibitem{Nicolai:1976}
H. Nicolai,
Supersymmetry and spin systems,
 \emph{J. Phys. A: Math. Gen.} \textbf{9}  (1976), 1497.

\bibitem{Rogers:2007}
A. Rogers,
\href{http://dx.doi.org/10.1142/1878}{Supermanifolds. Theory and applications},
\emph{Singapore: World Scientific},ISBN 978-981-02-1228-5/hbk, xii, 251 p. (2007).


\bibitem{Salam:1974}
A.~Salam \& J.~Strathdee, Super-gauge transformations,
\href{https://doi.org/10.1016/0550-3213(74)90537-9}{\emph{Nucl. Phys. B}} \textbf{76}  (1974), 477--482.

\bibitem{Schmitt:1997}
 T. Schmitt, Supergeometry and quantum field theory, or: what is a classical configuration?, \href{https://doi.org/10.1142/S0129055X97000348}{{\emph{Rev. Math.
Phys.}} } \textbf{9}, 993--1052, 1997.

\bibitem{Shvarts:1984}
A.S. Schwarz,
\newblock On the definition of superspace,
\newblock {\em Teoret. Mat. Fiz.}, \textbf{60}(1):37--42, 1984.

\bibitem{Varadrajan:2004}
 V.S.~Varadarajan, \href{http://bookstore.ams.org/cln-11}{Supersymmetry for mathematicians: an introduction},
Courant Lecture Notes in Mathematics, 11. New York University, Courant Institute of Mathematical Sciences, New York; \emph{American Mathematical Society}, Providence, RI, 2004. viii+300 pp. ISBN: 0-8218-3574-2.

  
  \bibitem{Volkov:1959}
 D. V. Volkov, 
 On the quantization of half-integer spin fields, \emph{Soviet Physics. JETP} \textbf{9} (1959), 1107--1111

\bibitem{Voronov:1984}
A.A. Voronov,
\newblock Maps of supermanifolds,
\newblock {\em Teoret. Mat. Fiz.}, \textbf{60}(1), 43--48, 1984.

\bibitem{Wess:1974}
J.~Wess \& B.~Zumino,
Supergauge transformations in four dimensions,
\href{https://doi.org/10.1016/0550-3213(74)90355-1}{\emph{Nuclear Phys.}} \textbf{B70} (1974), 39--50. 


\bibitem{Witten:1981}
E. Witten, 
Dynamical breaking of supersymmetry, \emph{Nucl. Phys. B} \textbf{188} (1981), 513.

\bibitem{Witten:2019}
E. Witten,
Notes On Supermanifolds and Integration,
\emph{Pure Appl.Math.Quart.} \textbf{15} (2019) 1, 3--56.
\end{small}
\end{thebibliography}
\end{document}